\documentclass{amsart}
\usepackage{amssymb,amsfonts,latexsym}

\newtheorem{theorem}{Theorem}[section]
\newtheorem{lemma}[theorem]{Lemma}
\newtheorem{proposition}[theorem]{Proposition}

\theoremstyle{definition}

\newcommand{\Ad}{\text{Ad}}

\newcommand{\gr}{\text{gr}}

\newcommand{\ot}{\otimes}

\newcommand{\ben}{\begin{enumerate}}
\newcommand{\een}{\end{enumerate}}

\newcommand{\Z}{{\mathbb Z}}

\newcommand{\C}{{\mathbb C}}

\begin{document}

\title[Quasi-Hopf algebras with radical of prime
codimension] {Liftings of graded Quasi-Hopf algebras with radical
of prime codimension}
\author{Pavel Etingof}
\address{Department of Mathematics, Massachusetts Institute of Technology,
Cambridge, MA 02139, USA} \email{etingof@math.mit.edu}

\author{Shlomo Gelaki}
\address{Department of Mathematics, Technion-Israel Institute of
Technology, Haifa 32000, Israel}
\email{gelaki@math.technion.ac.il}


\maketitle

\section{Introduction}

Let $p$ be a prime, and let $RG(p)$ denote the set of equivalence
classes of radically graded finite dimensional quasi-Hopf algebras
over $\mathbb C$, whose radical has codimension $p$. In
\cite{eg1},\cite{eg2} we completely describe the set $RG(p)$.
Namely, we show that for $p>2$, $RG(p)$ consists of the quasi-Hopf
algebras $A(q)$ constructed in \cite{g} for each primitive root of
unity $q$ of order $p^2$, the Andruskiewitsch-Schneider Hopf
algebras \cite{as1}, and semisimple quasi-Hopf algebras $H_+(p)$
and $H_-(p)$ of dimension $p$; on the other hand, $RG(2)$ consists
of the Nichols Hopf algebras of dimension $2^n$, $n\ge 1$, and the
special quasi-Hopf algebras $H(2)$, $H_\pm(8)=A(\pm i)$, and
$H(32)$ (of dimensions 2,8,32) constructed in \cite{eg1}.

The purpose of this paper is to classify finite dimensional quasi-Hopf
algebras $A$ whose radical is a quasi-Hopf ideal and
has codimension $p$; that is, $A$ with $\gr A\in RG(p)$,
where $\gr A$ is the associated graded algebra
taken with respect to the radical filtration on $A$.
The main result of this paper is the following theorem.

\begin{theorem}\label{main}
Let $A$ be a finite dimensional quasi-Hopf algebra whose radical
is a quasi-Hopf ideal of
prime codimension $p$. Then either $A$ is twist equivalent to a
Hopf algebra, or it is twist equivalent to $H(2)$, $H_\pm(p)$,
$A(q)$, or $H(32)$.
\end{theorem}

Note that if $A$ is a Hopf algebra (up to twist) then for $p=2$
$A$ is a Nichols Hopf algebra of dimension $2^n$, while for $p>2$
it is the dual to a finite dimensional pointed Hopf algebra with
group of grouplike elements $\Bbb Z_p$. Such Hopf algebras are
completely classified for $p\ne 5,7$ in \cite{as2} (see Remark
1.10(v), \cite{mu}).

Note also that any finite tensor category whose simple objects are
invertible and form a group of order $p$ under tensor is the
representation category of a quasi-Hopf algebra $A$ as above. Thus
this paper provides a classification of such categories.

The organization of the paper is as follows. Section 2 is devoted
to preliminaries. In Section 3 we introduce a construction
of semidirect products for quasi-Hopf algebras, which is used in
the sequel. In Section 4 we consider
liftings of the quasi-Hopf algebras $A(q)$ from [G],
and of $H(32)$. We prove:

\begin{theorem}\label{main1}
Let ${\rm gr}(A)$ be $A(q)$ or $H(32)$.
Then $A$ is equivalent to ${\rm gr}(A)$ as a quasi-Hopf algebra.
\end{theorem}

In Section 5 we consider liftings of Andruskiewitsch-Schneider
Hopf algebras \cite{as1}. We prove:

\begin{theorem}\label{main2}
Let ${\rm gr}(A)$ be an Andruskiewitsch-Schneider Hopf algebra.
Then $A$ is twist equivalent to a Hopf algebra.
\end{theorem}

By \cite{eg1}, \cite{eg2}, Theorem \ref{main} follows from
Theorems \ref{main1} -- \ref{main2}.

{\bf Acknowledgments.} The research of the first author was
partially supported by the NSF grant DMS-9988796.
The second author was supported by Technion V.P.R.
Fund - Dent Charitable Trust - Non Military Research Fund, and by The Israel
Science Foundation (grant No. 70/02-1). He also thanks MIT for its
warm hospitality.
Both authors were supported by BSF grant No. 2002040.

\section{Preliminaries}

All constructions in this paper are done over the field of complex
numbers $\mathbb{C}$.

We refer the reader to \cite{d} for the definition of a quasi-Hopf
algebra and a twist of a quasi-Hopf algebra.

\subsection{The quasi-Hopf algebras $A(q)$}
\begin{theorem} \cite{g}\label{propp3}
Let $n\ge 2$ be an integer.

\noindent
 There exist $n^3-$dimensional quasi-Hopf algebras $A(q)$, parametrized
by primitive roots of unity $q$ of order $n^2$,
which
have the following structure. As algebras $A(q)$ are generated
by $a,x$ with the relations $ax=q^n xa$, $a^n=1$, $x^{n^2}=0$.
The element $a$
is grouplike, while the coproduct of $x$ is given by the formula
\begin{equation*}
\Delta(x)=x\ot \sum_{y=0}^{n-1}q ^y {\bf 1}_y + 1\ot (1-{\bf
1}_{0})x + a^{-1}\ot {\bf 1}_{0}x,
\end{equation*}
where $\{\mathbf{1}_i|0\le
i\le n-1\}$ is the set of primitive idempotents of $\C[a]$
defined by the condition $a{\bf 1}_i=q ^{ni}{\bf 1}_i$. The associator
is
$$\Phi:=\sum_{i,j,k=0}^{n-1}q^{\frac{-i(j+k-(j+k)')}{n}}{\bf 1}_i\ot
{\bf 1}_j\ot {\bf 1}_k,$$ the distinguished
elements are $\alpha=a$,
$\beta=1$, and the antipode is $S(a)=a^{-1}$,
$S(x)=-x\sum_{z=0}^{n-1}q ^{n-z}{\bf 1}_z$.
\end{theorem}

Note that $S^2(x)=q^{n+1}x$, $S^2(a)=a$, so $S^{2n}=\Ad(a)$.

\subsection{The quasi-Hopf algebra $H(32)$}
\begin{theorem} \cite{eg1} \label{32}
There exists a $32-$ dimensional quasi-Hopf algebra $H(32)$, which
has the following structure. As an algebra $H(32)$ is generated by
$a,x,y$ with the relations $ax=-xa$, $ay=-ya$, $a^2=1$, $x^4=0$,
$y^4=0$, $xy+iyx=0$. The element $a$ is grouplike, while the
coproducts of $x,y$ are given by the formulas $$\Delta(x)=x\otimes
(p_+ + ip_-) +1\otimes p_+x + a\otimes p_-x,$$
$$\Delta(y)=y\otimes (p_+ - ip_-) +1\otimes p_+y + a\otimes
p_-y,$$ where $p_+:=(1+ a)/2, p_-:=(1-a)/2$. The associator is
$\Phi=1-2p_-\otimes p_-\otimes p_-$, the distinguished elements
are $\alpha=a$, $\beta=1$, and the antipode is $S(a)=a$,
$S(x)=-x(p_+ + ip_-)$, $S(y)=-y(p_+-ip_-)$. Thus, $H(32)$ is
generated by its quasi-Hopf subalgebras $A(i)$ and $A(-i)$
generated by $a,x$ and $a,y$, respectively.
\end{theorem}

Note that $S^2(x)=ix$ and $S^2(y)=-iy$, so $S^4=\Ad(a)$.

\subsection{Andruskiewitsch-Schneider Hopf algebras}

In \cite{as1}, Theorem 1.3, Andruskiewitsch and Schneider
completely classified finite dimensional coradically graded Hopf
algebras with coradical of prime dimension $p\ge 3$. Since this
collection of Hopf algebras is obviously self-dual, the list of
Theorem 1.3 of \cite{as1} also classifies radically graded Hopf
algebras with radical of prime codimension $p\ge 3$.

Let us recall this list.

1. Quantum lines (Taft Hopf algebras) of dimension $p^2$ (type
$A_1$). They are generated by $a,x$ where
$$
a^p=1,\ x^p=0,\ ax=qxa,
$$
and
$$
\Delta(a)=a\otimes a, \Delta(x)=x\otimes a+1\otimes x.
$$
Here $q$ is a primitive root of unity of degree $p$.

2. Quantum planes (book Hopf algebras) of dimension $p^3$
(type $A_1\times A_1$).
They are generated by $a,x,y$ where
$$
a^p=1,\ x^p=0,\ y^p=0,\ xy-yx=0,\ ax=qxa,\ ay=q^mya
$$
and
$$
\Delta(a)=a\otimes a, \Delta(x)=x\otimes a+1\otimes x,
\Delta(y)=y\otimes 1+a^m\otimes y.
$$
Here $q$ is primitive a root of unity of degree $p$,
and $m\in \Bbb Z_p$ is a nonzero element.

3. Algebras of type $A_2,B_2,G_2$. They are generated by $a,x,y$
satisfying the relations $axa^{-1}=qx$ and $aya^{-1}=q^by$, where
$b\in \Z_p$ is a nonzero element, and additional Serre relations
coming from the Nichols algebra ${\mathfrak u}_q^+$ (of type
$A_2,B_2$, or $G_2$, respectively) generated by $x,y$. The element
$a$ is grouplike, while the coproducts of $x,y$ are given by the
formulas
$$\Delta(x)=x\ot a+1\ot x,\;\Delta(y)=y\ot a^d+1\ot y,$$ where
$0\ne d\in \Z_p$, $b+d=-1$, and the following conditions hold:

1) The case $A_2$: $p=3$, or $p=1$ modulo 3, and $d^2+d+1=0$.

2) The case $B_2$: $p=1$ modulo 4, and $2d^2+2d+1=0$.

3) The case $G_2$: $p=1$ modulo 3, and $3d^2+3d+1=0$.

Here, as before, $q$ is a primitive root of unity of degree $p$.

4. Algebras of type $A_2\times A_1$ ($p=3$). They are generated by
$a,x,y$ as above (with $b=d=1$), and a new element $z$ such that
$aza^{-1}=q^fz$, $$\Delta(z)=z\otimes a^{-f}+1\otimes z,$$ where
$f=1$ or $2$ (subject to appropriate Serre relations).

5. Algebras of type $A_2\times A_2$ ($p=3$). They are generated by
$a,x,y$ as above (with $b=d=1$), and new elements $x',y'$, such
that $ax'a^{-1}=q^fx'$, $ay'a^{-1}=q^{f}y'$,
$$\Delta(x')=x'\otimes a^{-f}+1\otimes x',\;\Delta(y')=y'\otimes
a^{-f}+1\otimes y',$$ where $f=1$ or $2$ (subject to appropriate
Serre relations).

In all cases the Hopf algebra is isomorphic
$\C[\Z_p]\ltimes {\mathfrak u}_q^+$ as an {\em algebra}
(where ${\mathfrak u}_q^+$ is the corresponding Nichols algebra).

\subsection{Deformations}

Let $A$ be a finite dimensional quasi-Hopf algebra, in which the
radical is a quasi-Hopf ideal. Let $A_0$ be the associated graded
quasi-Hopf algebra of $A$ under the radical filtration. The
algebra $A_0$ has a grading $A_0=\bigoplus_{i=0}^{N} A_0[i]$ by
nonnegative integers. Thus we have gradings on the spaces of
linear maps between tensor powers of $A_0$. If $f_0$ is such a map
of some degree $d$ and $f$ is obtained from $f_0$ by adding terms
of degree higher than $d$, we will write $f=f_0+{\rm hdt}$ (where
hdt stands for higher degree terms) and say that $f$ is a
deformation (or lifting) of $f_0$. \footnote{More generally, if
$f_1,...,f_n$ have degrees $d_1,...,d_n$ then $f_1+...+f_n$ + hdt
will mean the sum of $f_i$ plus a sum of corrections of degrees
higher than all $d_i$.} Thus we can talk about deformations of
algebra, coalgebra structures, associator, antipode, etc. In
particular, we may say that $A$ is a deformation of $A_0$, in the
sense that $A$ is identified with $A_0$ as a graded vector space,
so that the structure maps of $A$ are given by some deformations
of the structure maps of $A_0$. This is the way we will think
about $A$ in the considerations below.

\subsection{Quasi-Hopf liftings of Hopf algebras}

\begin{proposition} \label{qhl}
Let $A$ be a finite dimensional quasi-Hopf algebra whose
radical is a quasi-Hopf ideal, such that $A_0={\rm gr}(A)$ is a
Hopf algebra (i.e., the associator of $A_0$ is equal to $1$). If
$H^3(A_0^*,\Bbb C)=0$ then $A$ is twist equivalent to a Hopf
algebra.
\end{proposition}

\begin{proof}
Let $\Phi$ be the associator of $A$.
We have $\Phi=1+$ hdt. Write $\Phi=1+\gamma +$ hdt, where $\gamma$
is of degree $d$. Then $\gamma$ is a $3-$cocycle of $A_0^*$ with
coefficients in the trivial $A_0^*-$bimodule $\C$. Since
$H^3(A_0^*,\Bbb C)=0$, we have $\gamma=dj$, where $j\in
A_0^{\otimes 2}[d]$, and the twisted associator $\Phi^{1+j}$ is
equal to $1+$ terms of degree $\ge d+1$. By continuing this
procedure, we will come to a situation where $\Phi=1$, as desired.
\end{proof}

\section{Semidirect products}

Let $H$ be a quasi-Hopf algebra with associator $\Phi$,
 $g:H\to H$ an algebra automorphism.
Define $\Delta_{g}(h):=(g\ot g)(\Delta(g^{-1}(h)))$.
Assume that
\begin{equation}\label{autg}
\Delta_{g}(h)=
K\Delta(h)K^{-1},\;h\in H,
\end{equation}
for some invertible element $K\in H\ot H$ such that
$(\varepsilon\otimes 1)(K)=(1\otimes \varepsilon)(K)=1$, and
$$
g^{\ot 3}
(\Phi)=\Phi^K:=(1\ot K)(id\ot \Delta)(K)\Phi(\Delta\ot
id)(K^{-1})(K^{-1}\ot 1).
$$

\begin{lemma} The semidirect product algebra $\C[g,g^{-1}]
\ltimes H$ is a quasi-Hopf algebra with associator $\Phi$,
$\Delta(g)=K^{-1}(g\ot g)$ and $\varepsilon(g)=1$.
\end{lemma}

\begin{proof}
Using equation (\ref{autg}), we see that
there exists a well defined coproduct on the algebra
$\C[g,g^{-1}]\ltimes H$ which extends the coproduct on $H$
and such that $\Delta(g)=K^{-1}(g\ot g)$.
It is straightforward to verify that $(\Delta\ot id)\Delta(g)=
\Phi^{-1}(id\ot \Delta)\Delta(g)\Phi$.
It is also easy to show that the antipode extends to the semidirect
product. The lemma is proved.
\end{proof}

Suppose further that $g^n(h)=aha^{-1}$ for some $a\in H$,
and that
\begin{equation}\label{jg}
(g^{n-1})^{\ot 2}(K)\dots g^{\ot 2}(K)K =(a\ot a)\Delta(a)^{-1}.
\end{equation}

\begin{lemma}\label{sd} The ideal generated by $g^n-a$ is a quasi-Hopf
ideal, so the quotient $\tilde H:=
(\C[g,g^{-1}]\ltimes H)/<g^n-a>$ is a quasi-Hopf algebra.
\end{lemma}

\begin{proof}
Using condition (\ref{jg}) it is straightforward to see that
\begin{eqnarray*}
\lefteqn{\Delta(g^n-a)}\\
& = &
\Delta(a)(a^{-1}g^n\ot a^{-1}g^n -1)
\\
& = & \Delta(a)(a^{-1}(g^n-a)\ot a^{-1}g^n +
1\ot a^{-1}(g^n-a)).
\end{eqnarray*}
It is also straightforward to verify that $<g^n-a>$ is invariant under the
antipode. Thus $<g^n-a>$ is a quasi-Hopf ideal.
\end{proof}

\section{Proof of Theorem \ref{main1}}

Let $A_0=A(q)$, where $q$ is a
primitive
root of unity of order $n^2$,  or $A_0=H(32)$ (in which case we
set $n=2$). Let $A$ be a lifting of $A_0$.
We apply the semidirect product construction to the
automorphism $g:=S^2$ of the algebra $A$.

Proposition 1.2 of \cite{d} provides a twist $K$ such that $\Delta_{g}(y)=
K\Delta(y)K^{-1}$, for all $y\in A$. Let $K_0$ be the degree zero
part of $K$. Then $K_0$ commutes with $\Delta_0(y)$ for any $y\in
A_0$, where $\Delta_0$ is the comultiplication of $A_0$.
This implies, using straightforward computations, that $K_0=1$,
and hence $K=1+$ hdt.

\begin{lemma}
The automorphism $S^{2n}$ of $A$ is inner.
\end{lemma}

\begin{proof}
Let $G:={\rm Aut}(A)/{\rm Inn}(A)$ be the group of outer algebra
automorphisms of $A$; it is a linear algebraic group. Let $S_0$ be
the antipode of $A_0$. Then $S_0^{2n}={\rm Ad}(a)$. Hence the
class $[S^{2n}]$ of $S^{2n}$ is unipotent in $G$. On the other
hand, from the main result of \cite{eno} (see also \cite{hn}) we
know that on the category of representations of $A$ there exists
an isomorphism of tensor functors $\theta_V:V\to \chi\ot V^{****}\ot
\chi^{-1}$, where $\chi$ is a $1-$dimensional representation of
$A$. In our case, $\chi$ has order $n$. Therefore there exists a
tensor isomorphism $c_V:V\to V^{**^{2n}}$, and thus $S^{4n}$ is an
inner automorphism of $A$: $S^{4n}(z)=czc^{-1}$ where
$c\in A$ is the element realizing the automorphism $c_V$ (i.e., $c_V=c|_V$).
Hence, $[S^{2n}]=1$ in $G$ and $S^{2n}$
is inner.
\end{proof}

Thus, $S^{2n}(z)=b zb^{-1}$ for some element $b\in
A$, where $b=a+$ hdt is a deformation of the grouplike
element $a\in A_0$.

\begin{lemma}\label{choi} The element $b$ can be chosen so that
\begin{equation}\label{jg1}
K^{-1}g^{\ot 2}(K^{-1})\dots (g^{n-1})^{\ot
2}(K^{-1})=\Delta(b)(b^{-1}\ot b^{-1}).
\end{equation}
\end{lemma}

\begin{proof} Denote the left hand side of (\ref{jg1}) by $T$.

Let $c$ be the element defined above such that $S^{4n}(z)=czc^{-1}$.
It is easy to check that $c=a^2+$ hdt.

We claim that
$$
T(b\otimes b)T(b^{-1}\otimes b^{-1})=\Delta(c)(c^{-1}\otimes c^{-1}).
$$
Indeed, since $\{c_V\}$ is an isomorphism of tensor functors, $L(c_V\ot c_W)=
c_{V\ot W}$, where $L$ is the tensor structure on the functor
$**^{2n}$. This defines an element $L\in A\ot A$ such that $L=\Delta(c)
(c^{-1}\ot c^{-1})$. However, we know from \cite{d} that the tensor structure
on $**$ is given by $K^{-1}$, and hence the tensor structure on $**^{2n}$
is given by $K^{-1}(S^2\ot S^2)(K^{-1})\cdots (S^{4n-2}\ot S^{4n-2})(K^{-1})$,
which is equal to $T(S^{2n}\ot S^{2n})(T)=T\Ad(b\ot b)(T)$, as desired.

Thus we have that $z:=b^2c^{-1}=1+$ hdt is a central element.
Replacing $b$ with $bz^{-1/2}$, we may assume without
loss of generality that $b$ is chosen to satisfy the condition
$z=1$, i.e. $b^2=c$.

We claim that for this choice of $b$,
equation (\ref{jg1}) holds.
Indeed, assume the contrary,
and let the lowest degree of nontrivial terms in the difference
between the two sides of
(\ref{jg1}) be $d$. Let $\gamma\ne 0$ be the part of the difference that is
homogeneous of degree $d$; so $T=\Delta(b)(b^{-1}\otimes
b^{-1})-\gamma$ modulo elements of degree $\ge d+1$.
Then modulo elements of degree $\ge
d+1$, we have for any $y\in A$,
\begin{eqnarray*}
\lefteqn{\Delta_{{\rm Ad}(b)}(y)}\\ & = &
(b\otimes b)\Delta(b^{-1}yb)(b^{-1}\otimes b^{-1})=T^{-1}\Delta(y)T\\
& = & ((b\otimes b)\Delta(b)^{-1}+\gamma)
\Delta(y)(\Delta(b)(b^{-1}\otimes b^{-1})-\gamma).
\end{eqnarray*}
This means that for any $y\in A_0$ one has
$[\gamma,\Delta_0(y)]=0$, and in particular
$[\gamma,a\otimes a]=0$.

Also, modulo terms of degree $\ge d+1$,
$$
(\Delta(b)(b^{-1}\otimes b^{-1})-\gamma)(b\otimes
b)(\Delta(b)(b^{-1}\otimes b^{-1})-\gamma)
(b^{-1}\otimes b^{-1})=\Delta(c)(c^{-1}\otimes c^{-1}).
$$
Since $b^2=c$, this implies that
$$
\gamma+(a\otimes a)\gamma(a^{-1}\otimes a^{-1})=0.
$$
Since $[\gamma,a\otimes a]=0$, this is equivalent to
$2\gamma=0$. Contradiction.
\end{proof}

Lemma \ref{sd} and Lemma \ref{choi} imply that we can define the
semidirect product quasi-Hopf algebra $$\tilde
H:=(\C[g,g^{-1}]\ltimes A)/<g^n -b>;$$ it is of dimension $n^4$
for $A=A(q)$, and 64 for $A=H(32)$.

Let $J\in \tilde H\ot \tilde H$ be a deformation of the twist
$J_0$ from \cite{g}. Then it follows from \cite{g} that
$H:={\tilde H}^ {J^{-1}}$ is a quasi-Hopf lifting (=deformation)
of a Hopf algebra $H_0$. If $A_0=A(q)$ then $H_0$ is the Taft Hopf
algebra of dimension $n^4$. On the other hand, if $A_0=H(32)$ then
$H_0$ is a book Hopf algebra of dimension $64$, generated by
$g,x_+,x_-$ satisfying $gx_+g^{-1}=ix_+$, $gx_-g^{-1}=-ix_-$,
$x_+x_-=x_-x_+$.  The element $g$ is grouplike and
$$\Delta(x_+)=x_+\ot g+1\ot x_+,\;\Delta(x_-)=x_-\ot 1+g^{-1}\ot x_-.$$

\begin{theorem}\label{taft}
If $H$ is a basic quasi-Hopf algebra such that $H_0=\gr H$
then there exists a quasi-Hopf twist $F$
such that $H\cong H_0^F$ and $F=1$ modulo higher terms.
\end{theorem}

\begin{proof}
Let us first show that $H$ is twist equivalent to a Hopf algebra.
By Proposition \ref{qhl}, for this it is sufficient to show that
$H^3(H_0^*,\Bbb C)=0$. If $H_0$ is the Taft algebra then
$H_0^*\cong H_0$ and one finds (see  \cite{gk}), that
$H^\bullet(H_0^*,\C)$ is the polynomial algebra on an element of
degree 2. In particular, $H^3(H_0,\C)=0$, as required. On the
other hand, if $H_0$ is the book algebra of dimension $64$ then
$H_0^*$ is generated by $g,x,y$ with $g^4=1$, $x^4=y^4=xy-yx=0$,
$gxg^{-1}=ix$, $gyg^{-1}=iy$. From this it is easy to obtain using
K\"unneth formula that $H^\bullet(H_0^*,\Bbb C)$ is the polynomial
algebra in two elements of degree 2 (see \cite{gk}). In
particular, again $H^3(H_0^*,\Bbb C)=0$.

Now, it is easy to see, using the
methods of \cite{as2}, that in both cases
(the Taft algebra and the book algebra),
$H_0$ has no non-trivial Hopf liftings,
so there exists a twist $F=1+$ hdt and an isomorphism
of quasi-Hopf algebras $\eta:H_0^{F}\to H$, as desired.
\footnote{We note that if $H_0$ is the book algebra of
dimension 64, then $H_0^*$, unlike $H_0$,
has a nontrivial Hopf lifting, which
is the dual of the quantum $sl_2$ at the 4-th root of unity.}
\end{proof}

Replacing $J$ with $JF$, we may assume that $F=1$. Thus we find
that $\tilde H=H_0^{J}$. So we have a homomorphism of quasi-Hopf
algebras
\begin{equation}\label{eta}
\eta: A\to H_0^{J}.
\end{equation}

The map $\eta$ may be regarded as an injective linear map $A_0\to
H_0^J$, whose image is a quasi-Hopf subalgebra of $H_0^J$, and
whose degree $0$ part is the ``tautological'' homomorphism
$\eta_0: A_0\to H_0^{J_0}$\footnote{Note the difference between
the map $\eta_0: A_0\to H_0$ and the map $A_0\to H_0^*$ (not
$H_0$!) constructed in \cite{eg1}.}. Let $E$ be the set of pairs
$(\eta,J)$ satisfying this condition. On $E$ there is an action of
the ``gauge group'' ${\mathcal G}$, generated by the following
three types of transformations.

1. $(\eta,J)\to (\eta\circ \xi,J)$, where $\xi=1$+hdt is a linear
automorphism of $A_0$.

2. $(\eta,J)\to ({\rm Ad}(h)\circ \eta,(h\otimes
h)J\Delta(h)^{-1})$,
where $h=1+$hdt is an element of $H_0$, and $\Delta$ is the
coproduct of $H_0$.

3. $(\eta,J)\to (\eta,TJ)$, where $T=1$+hdt belongs to
$\eta(A_0)^{\otimes 2}$.

It is clear that the transformations from ${\mathcal G}$ do not
change the equivalence class of the quasi-Hopf algebra
$\eta(A_0)$. Thus, Theorem \ref{main1}
follows from the following proposition.

\begin{proposition}
Any pair $(\eta,J)\in E$ is equivalent under
${\mathcal G}$ to $(\eta_0,J_0)$.
\end{proposition}

\begin{proof} It suffices to prove that for any $d>0$,
any pair $(\eta,J)\in E$ can be brought by the action of
${\mathcal G}$ to a pair equal to $(\eta_0,J_0)$ modulo terms of
degree $\ge d$. We will prove it by induction in $d$. The base of
induction is obvious. To prove the induction step, we need to take
a pair $(\eta,J)$ which is equal to $(\eta_0,J_0)$ modulo terms of
degree $d$, and show it is equivalent to a pair equal to
$(\eta_0,J_0)$ modulo terms of degree $d+1$.

Thus, the multiplication in $A_0$ obtained by pullback of the
multiplication of $H_0$ by $\eta$
is given by $z*w=zw+c_d(z,w)$+hdt, where $zw$ is the
multiplication in $A_0$, and $c_d$ is a bilinear operator
$A_0\times A_0\to A_0$ of degree $d$, which is a Hochschild 2-cocycle.
Similarly, $\eta(z)=\eta_0(z)+\eta_d(z)+$ hdt, where $\eta_d$ is a degree $d$
operator $A_0\to H_0$, and $J=TJ_0$, $T=1+t_d+$ hdt, where
$t_d\in H_0\otimes H_0$ has degree $d$.

We have $\eta(z*w)=\eta(z)\eta(w)$ which is equivalent to
$$(\eta_0+\eta_d+\cdots)(zw+c_d(z,w)+\cdots)=(\eta_0+\eta_d+\cdots)
(z)(\eta_0+\eta_d+\cdots)(w).$$ Taking $z$ of degree $d_z$ and $w$
of degree $d_w$, and looking at terms of degree $d+d_z+d_w$ we get
$\eta_d(zw)+\eta_0(c_d(z,w))=\eta_d(z)\eta_0(w)+\eta_0(z)\eta_d(w)$,
which implies $\eta_0(c_d(z,w))=\eta_d(z)\eta_0(w)+\eta_0(z)\eta_d(w)-\eta_d(zw)$.
Therefore, $\eta_0(c_d)$ is the differential of a 1-cochain $\eta_d\in
C^1(A_0,H_0)$. So $\eta_0$ sends the
cohomology class of $c_d\in H^2(A_0,A_0)$ into 0 in
$H^2(A_0,H_0)$. Since $A_0$ is a direct summand in $H_0$ as an
$A_0$-bimodule, this implies that the class $[c_d]\in H^2(A_0,A_0)$ is
zero. Thus, $c_d$ is the differential of some element $\tau$ of degree
$d$ in $C^1(A_0,A_0)$. So, using an element $g_1\in {\mathcal G}$ of
type $1$, we can replace $(\eta,J)$ with an equivalent pair for
which $c_d=0$. Thus in the sequel we may assume that $c_d=0$.

Now, $\eta_d$ is a Hochschild 1-cocycle of $A_0$ with coefficients
in $H_0$. Consider the cohomology class $[\eta_d]$ in
$H^1(A_0,H_0)$. We claim that in fact, $[\eta_d]$ belongs
to the subspace $H^1(A_0,A_0)\subset H^1(A_0,H_0)$. To prove
this, note that since $\eta(A_0)$ is closed under coproduct,
for any $z\in A_0$ we have
$$
\Delta_0(\eta_d(z))-(\eta_d\otimes \eta_0+\eta_0\otimes
\eta_d)(\Delta_0(z))=[\Delta_0(z),t]\; {\rm mod }\;A_0\otimes A_0,
$$
where $\Delta_0$ denotes the coproduct of $H_0^{J_0}$ and
$A_0$, and $t:=t_d$.

Now recall that $H_0=\oplus_{k=0}^{n-1}A_0g^k$.
Thus, $\eta_d=\sum_{k=0}^{n-1} \omega_k$, where
$\omega_k\in Z^1(A_0,A_0g^k)$ (a 1-cocycle). Similarly,
$t=\sum_{k,l=0}^{n-1} t_{kl}$, $t_{kl}\in A_0g^k\otimes A_0g^l$.
Thus for $k\ne 0$ we have
$$
-(\omega_k\otimes \eta_0)(\Delta_0(z))=[\Delta_0(z),t_{k0}]
$$
Applying the counit in the second component, we get that
$\omega_k(z)=[z,y_k]$
for some $y_k\in A_0g^k$, which implies
that $\omega_k$ is an inner derivation.

Thus $[\eta_d]\in
H^1(A_0,A_0)$, and hence we could have chosen $g_1$ so that
$[\eta_d]=0$. So we may assume $[\eta_d]=0$ and
thus $\eta_d$ is the bracketing operator with an element $y\in H_0$.
Thus using an appropriate element $g_2\in {\mathcal G}$ of type
2, we can assume in the sequel that $\eta_d=0$.

Now we will show that in the situation when $c_d=0$ and
$\eta_d=0$, we must have $t\in A_0\otimes A_0$. Then by using an
element $g_3\in {\mathcal G}$ of type 3, we can come to a
situation with $t=0$ which completes the induction step.

Since the associator $\Phi:=\Phi_0^T$ lies in
$\eta(A_0)^{\otimes 3}$, we find that
\begin{equation}\label{dt}
t\ot 1 + (\Delta_0\ot id)(t)-\Phi_0 (id\ot \Delta_0)(t)
\Phi_0^{-1}+\Phi_0 (1\ot t) \Phi_0^{-1}\in A_0^{\otimes 3}.
\end{equation}
Let $H_0=A_0\oplus H'$, where $H'$ is the direct sum of $A_0g^k$ for
$k=1,...,n-1$.
Let $t'$ be the projection
of $t$ to $H'\otimes H_0$. Since the first component of
$t-t'$ is in $A_0$, we get from (\ref{dt}):
\begin{equation}\label{dt'}
t'\ot 1 + (\Delta_0\ot id)(t')-\Phi_0 (id\ot \Delta_0)(t')
\Phi_0^{-1}+\Phi_0 (1\ot t') \Phi_0^{-1}\in A_0\otimes H_0\otimes H_0.
\end{equation}
But out of the four summands in this expression,
all but the last one belong to $H'\otimes H_0\otimes H_0$, while the
last summand is in $A_0\otimes H_0\otimes H_0$.
Hence the last summand is zero and $t'=0$.
Similarly, one shows that the projection $t''$ of $t$ to
$H_0\otimes H'$ is zero. Hence, $t\in A_0\otimes A_0$, as
desired.
\end{proof}

The theorem is proved.

\section{Proof of theorem \ref{main2}}

Let $H$ be a Hopf algebra
from the list of Subsection 2.3.
Since this list is self-dual, Theorem \ref{main2}
follows from Proposition \ref{qhl} and the following proposition.

\begin{proposition} We have $H^3(H,\C)=0$.
\end{proposition}

\begin{proof}
By a standard argument, $H^3(H,\C)=H^3(\C[\Z_p]\ltimes {\mathfrak u}_q^+,\C)=
H^3({\mathfrak u}_q^+,\C)^{\Z_p}$.

Let us first assume that $p>3$. In this case, according to Section
2.3, the rank of the corresponding simple Lie algebra is $\le 2$,
and $p$ is greater than its Coxeter number. Thus we can apply
\cite{gk}, Theorem 2.5. In this theorem, it is shown that
$H^\bullet({\mathfrak u}_q^+,\C) =\bigoplus_{w\in W}\C {\eta_w}\ot
S({\frak n}_+)$, where $S({\frak n}_+)$ is the symmetric algebra
of the positive part of the associated Lie algebra, sitting in
degree 2, $W$ is the Weyl group, and $\eta_{w}$ has degree $l(w)$
(=length of $w$) .

The action of the generator
$a$ of $\Z_p$ on $H^\bullet({\mathfrak u}_q^+,\C)$ is given as follows: $a$ acts
trivially on $S({\frak n}_+)$ and $a\circ \eta_{w}=\lambda_w\eta_w$, where
in the $A_1$-case $\lambda_w=q^{-1}$ for the nontrivial element
$w\in W$, and in the other (rank 2) cases
$\lambda_w=q^{-(m+nd)}$, where $m,n$ are determined by writing
the weight $\gamma_w:=
\sum_{\alpha\in R^+,\;w(\alpha)<0}\alpha=\rho-w(\rho)$ as an
integral combination $m\alpha_1+n\alpha_2$ of the simple roots (here
$\rho$ is half the sum of positive roots, and $R_+$ is the set of
positive roots; in the non-simply laced case,
$\alpha_2$ always denotes the short root).

We claim that $H^3({\mathfrak u}_q^+,\C)^{\Z_p}=0$. Indeed, since $\Z_p$ acts
trivially on $S({\frak n}_+)$, and since $\lambda_w\ne 1$ if $w$
is a simple reflection, we have $H^3({\mathfrak u}_q^+,\C)^{\Z_p}=
(\bigoplus_{w\in W,\;l(w)=3}\C\eta_w)^{\Z_p}$. We will now show that
in each of the four cases $A_1\times A_1$, $A_2$, $B_2$ and $G_2$,
$\lambda_w\ne 1$, and hence $H^3({\mathfrak u}_q^+,\C)^{\Z_p}=0$.

In the case $A_1\times A_1$ there are no elements $w$ of length 3, so the result
is clear.

In the case $A_2$, there is one element $w$ of length 3 and one computes
that $\gamma_w=2\alpha_1+2\alpha_2$. Therefore $\lambda_w=(q^{-1})^{2+2d}$.
Since $d^2+d+1=0$, $d+1\ne 0$, so $\lambda_w\ne 1$.

In the case $B_2$, there are two elements $w_1,w_2$ of length 3 and one computes
that $\gamma_{w_1}=3\alpha_1+3\alpha_2$, $\gamma_{w_2}=2\alpha_1+4\alpha_2$.
Therefore $\lambda_{w_1}=(q^{-1})^{3+3d}$ and $\lambda_{w_2}=(q^{-1})^{2+4d}$.
Since $2d^2+2d+1=0$, we have $d+1,2d+1\ne 0$, so $\lambda_{w_1},\lambda_{w_2}\ne 1$.

In the case $G_2$, there are two elements $w_1,w_2$ of length 3 and one computes
that $\gamma_{w_1}=4\alpha_1+4\alpha_2$, $\gamma_{w_2}=2\alpha_1+6\alpha_2$.
Therefore $\lambda_{w_1}=(q^{-1})^{4+4d}$ and $\lambda_{w_2}=(q^{-1})^{2+6d}$.
Since $3d^2+3d+1=0$, we have $d+1,3d+1\ne 0$, so $\lambda_{w_1},\lambda_{w_2}\ne 1$.
Thus the result for $p>3$ follows.

Finally, let us consider the case $p=3$. In the rank 1 case the
above argument applies, so we consider rank $\ge 2$. Then the
argument above does not quite apply, since the rank of the Lie
algebra can now be bigger than 2, and $p$ may be equal to the
Coxeter number (a case not covered by Theorem 2.5 in \cite{gk}).
Thus we will use a slightly different argument. By \cite{gk}, the
Nichols algebra ${\frak u}_q^+$ has a filtration under which
$H^\bullet (\gr {\mathfrak u}_q^+,\C)=\Lambda_q\ot S({\frak
n}_+)$, where $\Lambda_q$ is the $q-$analogue of the exterior
algebra generated by $\varepsilon_{\alpha}$, $\alpha$ a positive
root, such that
$\varepsilon_{\alpha}\varepsilon_{\beta}+q^{<\alpha,\beta>}
\varepsilon_{\beta}\varepsilon_{\alpha}=0$ for $\beta<\alpha$ and
$\varepsilon_{\alpha}^2=0$ for any $\alpha$.

The generator $a$ of $\Z_p$ acts trivially on $S({\frak n}_+)$ and
by $a\varepsilon_{\alpha} a^{-1}=q^{-(m+nd)}\varepsilon_{\alpha}$,
where $\alpha=m\alpha_1+n\alpha_2$. The group $G$ of grouplikes in
the Frobenius-Lusztig quantum group acts on both $H^\bullet(\gr
{\mathfrak u}_q^+,\C)$ and $H^\bullet({\mathfrak u}_q^+,\C)$, and
the second one is isomorphic to a subrepresentation of the first
one by a spectral sequence argument \cite{gk}.

In the case $A_2$ ($d=1$) we claim that $H^3(\gr {\mathfrak u}_q^+,\C)^{\Z_3}=0$.
Indeed, we have
$H^3(\gr {\mathfrak u}_q^+,\C)=<\varepsilon_{\alpha_1}\varepsilon_{\alpha_2}
\varepsilon_{\alpha_1+\alpha_2}>$, and the generator $a$ of $\Z_p$ acts on this $1-$
dimensional space by $q^{-(1+d+1+d)}=q^{-4}\ne 1$.

In the case $A_2\times A_1$ ($d=1$), we have that the Nichols
algebra
${\frak u}_q^+$ is generated
by $x,y,z$ with $aza^{-1}=q^fz$. Thus, if $\beta$ is the simple
root corresponding to $z$, we have,
$H^3(\gr {\mathfrak u}_q^+,\C)^{\Z_3}
=<\varepsilon_{\beta}\varepsilon_{\alpha_1}\varepsilon_{\alpha_2}>$ if $f=2$, and is
zero if $f=1$.

In the case $A_2\times A_2$ ($d=1$), we have that the Nichols
algebra ${\frak u}_q^+$ is generated by $x,y,x',y'$ with
$axa^{-1}=qx$, $aya^{-1}=qy$, $ax'a^{-1}=q^fx'$ and
$ay'a^{-1}=q^fy'$. Hence, if $\beta_1,\beta_2$ are the roots
attached to $x',y'$, we find that $$H^3(\gr {\mathfrak
u}_q^+,\C)^{\Z_3}
=<\varepsilon_{\beta_i}\varepsilon_{\alpha_1}\varepsilon_{\alpha_2},
\varepsilon_{\alpha_i}\varepsilon_{\beta_1}\varepsilon_{\beta_2}|i=1,2>$$
if $f=1$, and is zero if $f=2$.

But even in cases where $H^3(\gr {\mathfrak u}_q^+,\C)^{\Z_3}\ne
0$, this cohomology is killed in the spectral sequence $H^\bullet(\gr
{\mathfrak u}_q^+,\C)^{\Z_3}\to H^\bullet({\mathfrak u}_q^+,\C)^{\Z_3}$
since by \cite{gk}, the only weights of terms under $G$ which can
survive in this sequence must be of the form $w(\rho)-\rho$ while
$-(\alpha_1+\alpha_2)\ne w(\rho)-\rho$ modulo 3. So in all cases
$H^3({\mathfrak u}_q^+,\C)^{\Z_3}=0$, and the result follows.
\end{proof}

\end{document}